\documentclass{gtpart}     
\usepackage{latexsym}
\usepackage{amsfonts}
\usepackage{amsbsy}
\usepackage{mathrsfs}
\usepackage{epsfig}
\usepackage[all]{xy}
\usepackage{url}

\title[Associahedra and Monoidal Categories]{Associahedra and Weak Monoidal Structures on Categories}

\author[Fiedorowicz]{Z. Fiedorowicz}
\givenname{Zbigniew}
\surname{Fiedorowicz}
\address{Department of Mathematics, The Ohio State University\\ Columbus, OH 43210-1174, USA}
\email{fiedorow@math.ohio-state.edu}
\urladdr{http://www.math.ohio-state.edu/people/fiedorow/view}

\author[Gubkin]{S. Gubkin}
\givenname{Steven}
\surname{Gubkin}
\address{Department of Mathematics, The Ohio State University\\ Columbus, OH 43210-1174, USA}
\email{gubkin@math.ohio-state.edu}
\urladdr{http://www.math.ohio-state.edu/people/gubkin/view}

\author[Vogt]{R.M.~Vogt}
\givenname{Rainer}
\surname{Vogt}
\address{Universit\"at Osnabr\"uck, Fachbereich Mathematik/Informatik\\ Albrechtstr. 28a, 49069 Osnabr\"uck, Germany}
\email{rainer@mathematik.uni-osnabrueck.de}
\urladdr{http://www.mathematik.uni-osnabrueck.de/staff/phpages/vogtr.rdf.shtml}

\keyword{monoidal categories}
\subject{primary}{msc2010}{18D10}
\keyword{operads}
\subject{secondary}{msc2010}{55P48}
\keyword{Tamari lattice}
\subject{secondary}{msc2010}{06A07}

\volumenumber{}
\issuenumber{}
\publicationyear{}
\papernumber{}
\startpage{}
\endpage{}
\doi{}
\MR{}
\Zbl{}
\received{}
\revised{}
\accepted{}
\published{}
\publishedonline{}
\proposed{}
\seconded{}
\corresponding{}
\editor{}
\version{}

\makeatletter
\@addtoreset{equation}{section}

\makeatother
\makeatletter
\@addtoreset{equation}{section}

\makeatother

\newtheorem{propo}{Proposition}[section]
\newtheorem{theo}[propo]{Theorem}

\newtheorem{coro}[propo]{Corollary}
\newtheorem{rema}[propo]{Remark}
\newtheorem{defi}[propo]{Definition}

\newtheorem{exam}[propo]{Example}
\newtheorem{lem}[propo]{Lemma}
\def\scK{\mathfrak{K}}
\def\scC{\mathcal{C}}
\def\scD{\mathcal{D}}
\def\scM{\mathcal{M}}
\def\scL{\mathcal{L}}
\def\scI{\mathcal{I}}
\def\hK{\widehat{\mathfrak K}}
\def\scO{\mathcal{O}}

\begin{document}

\begin{abstract}    
This paper answers the following question: what algebraic structure on a category corresponds to
an $A_n$ structure (in the sense of Stasheff) on the geometric realization of its nerve?
\end{abstract}

\maketitle

In his trailblazing paper \cite{Stash1}, Stasheff constructed an infinite hierarchy of higher homotopy associativity
conditions for an  H-space $X$. These conditions are parametrized by a family $\{K_n\}_{n\ge 2}$ of polyhedra,
which came to be known as associahedra.  The vertices of $K_n$ are in 1-1 correspondence with all possible
ways of associating an $n$-fold product $x_1x_2\dots x_n$, and an H-space $X$ is said to be an $A_n$-space
if there is a map $K_n\times X^n\longrightarrow X$ whose restriction to the vertices enumerates all possible
ways of associating the binary multiplication on $X$ into an $n$-fold multiplication. An $A_\infty$-space is known
to be equivalent to a strict monoid $MX$ and hence, up to group completion, to a loop space.

At the same time\footnote{Stasheff informs us that, although \cite{Stash1} and \cite{Mac2} both appeared in 1963,
Mac Lane's work preceded his and influenced his thinking. He further informs us (cf. \cite{Stash3}) that the associahedra were
implicitly defined in the even earlier work of Tamari \cite{Tam1}, \cite{Tam2}.}
 Mac Lane \cite{Mac2} analyzed higher associativity conditions for monoidal structures on categories.
He formulated analogs of Stasheff's $A_n$ conditions for categories. For $n=2,3,4$ the analogy is perfect.
In particular Mac Lane's $A_4$ condition is that a pentagonal diagram commute, whereas Stasheff's $K_4$ is a
pentagon. However for $n\ge 5$ the analogy breaks down.  Mac Lane's coherence theorem states that the $A_4$
condition implies all the higher $A_n$ conditions for $n\ge 5$. By contrast for any $n\ge 2$ one can construct
H-spaces $X$ which satisfy the $A_n$ condition but not the $A_{n+1}$ (or any higher) condition.

In this paper we show how Mac Lane's notion of a monoidal structure on a category can be weakened so as to
obtain a full hierarchy of $A_n$ conditions. The paper is similar in spirit to \cite{BFSV} where an $E_n$ hierarchy
of commutativity conditions on categories was considered, analogous to those on $n$-fold loop spaces. Similarly
to the case of associativity for categories, Joyal and Street \cite[Proposition 5.4]{JS} showed that if
these commutativity conditions are required to hold up to natural isomorphisms, then the $E_3$ condition
implies all higher $E_n$ conditions. In \cite{BFSV} we demonstrated that we could recover the entire $E_n$ hierarchy
for categories by weakening these commutativity conditions to hold up to natural transformations instead.
This strategy does not work for associativity, since LaPlaza \cite[Theorem 5]{LaPl} showed that even if the
associativity conditions are weakened to hold up to natural transformations, instead of isomorphisms,
this laxened form of Mac Lane's $A_4$ condition still implies all higher $A_n$ conditions. Thus a different
strategy for weakening Mac Lane's $A_n$ conditions for categories is required. 

In sections 1 and 2 we develop this strategy: we define the category theoretical ana\-logues of
Stasheff's associahedra in section 1 and $A_n$-monoidal categories in section 2. In section 3 we relate our work to
that of LaPlaza (and implicitly to that of Tamari) and give a simpler proof of his coherence result. In section 4 we prove a rectification result for
$A_\infty$-monoidal categories, similar in spirit to Mac Lane's rectification of a monoidal category to a strictly
monoidal one, by translating the rectification of an $A_\infty$-space to a monoid into category theory.

This paper presupposes some familiarity with the notion of operad and related concepts. A pr\'ecis of the
relevant definitions may be found in \cite{May} and some historical context in \cite{Stash2}. Since we will
be dealing exclusively with noncommutative operations, we will be using the non-$\Sigma$ forms of operads
throughout.

To forestall any possible misunderstanding, it should be pointed out that this paper is not related in any
significant way to the notion of $A_\infty$-category as developed by Fukaya, Kontsevich, Soibelman and others
(c.f. \cite{Keller} for an overview). In this paper we are discussing ordinary categories with weak monoidal
structures, not some notion of a weak higher category.

We would like to take this opportunity to thank Jim Stasheff and Stefan Forcey for some helpful suggestions
and references to previous work in this area.

\section{The associahedra as an operad in $CAT$}

 In order to keep track of  associativity data for our weakly monoidal categories, we will need a categorical 
equivalent of the associahedron $K_m$.  To begin with we formalize the notion of a parenthesized  word:

\begin{defi}
A \textit{parenthesized word} $(W,P)$is a finite linear order $W$ together with a (possibly empty) collection 
of closed intervals $P = \{p_i = [a_i, b_i]\}$ subject to the following requirements.
\begin{itemize}
\item The cardinality of each $p_i$ is at least 2 and is strictly smaller than the cardinality of $W$.
\item For any $i,j$, either $p_i \subset p_j$, $p_j \subset p_i$ or $p_i \cap p_j = \varnothing.$
\end{itemize}
\end{defi}

A parenthesized word $(W,P)$ can be converted into a parenthesized string of characters by putting as many left 
parentheses in front of an element $a\in W$ as $a$ is an initial element
of some $p_i\in P$ and as many right parentheses after an element $b\in W$ as $b$ is a final element in some $p_i\in P$, 
and concatenating the resulting characters.  For instance
$$\left\{x_1<x_2<x_3<x_4<x_5<x_6,\{[x_2,x_6],[x_2,x_4],[x_5,x_6]\}\right\}\mapsto x_1((x_2x_3x_4)(x_5x_6)).$$
It is clear that $(W,P)$ can be recovered from the parenthesized string and we will often find it convenient to 
represent $(W,P)$ in this way.
In most cases we will use the standard linear orders $W_m=\{x_1<x_2\dots <x_m\}$. In some induction arguments 
however we will need to consider subintervals of the $W_m$.

\begin{defi}\label{Kmdefi}
We define $\scK_m$ to be the poset of parenthesized words on the linear order $W_m$, where $(W_m,P_2) 
\leq (W_m, P_1)$ iff $P_1 \subset P_2$.  The minimal elements in this order are called the 
\textit{fully parenthesized} words of length $m$. In the degenerate cases $m=1$ and $m=0$, the poset $\scK_1$ 
consists of the single parenthesized word $id=(W_1,\emptyset)$, and $\scK_0$ consists of the single parenthesized
word $0=(\emptyset,\emptyset)$. The string $x_1x_2...x_m=(W_m,\emptyset)$ is the terminal object in $\scK_m$. 
As noted above, sometimes it will be convenient to use some other linear order
$W'$ of the same cardinality $m$. In that case the unique order isomorphism between $W'$ and $W_m$ specifies a 
canonical isomorphism between $\scK_m$ and the corresponding poset of parenthesized
words on $W'$.
\end{defi}
\pagebreak

\begin{exam}
The poset $\scK_4$:

$$
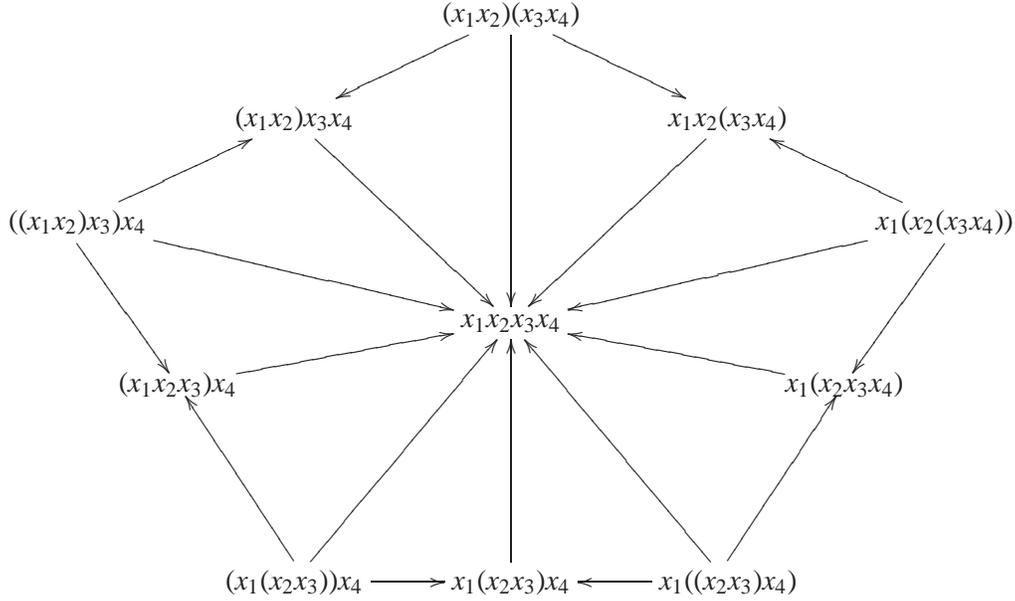
\xymatrix{
&&(x_1x_2)(x_3x_4)\ar[dl]\ar[dr]\ar[ddd]\\
&(x_1x_2)x_3x_4\ar[ddr] &&x_1x_2(x_3x_4)\ar[ddl]\\
{((x_1x_2)x_3)x_4}\save[]+<0pt,-8pt>*{}="a"\restore\ar[ur]\ar[drr] &&
 &&x_1(x_2(x_3x_4))\save[]+<0pt,-8pt>*{}="f"\restore\ar[ul]\ar[dll]\\
&&x_1x_2x_3x_4\\
{}\save[]+<38pt,10pt>*{(x_1x_2x_3)x_4}="b"\ar[urr]\restore &&
 &&{}\save[]+<-38pt,10pt>*{x_1(x_2x_3x_4)}="d"\ar[ull]\restore\\
\\
&\save[]+<0pt,8pt>*{}="c"\restore(x_1(x_2x_3))x_4\ar[r]\ar[uuur]
 &x_1(x_2x_3)x_4\ar[uuu] &\save[]+<0pt,8pt>*{}="e"\restore x_1((x_2x_3)x_4)\ar[l]\ar[uuul]
\ar "a";"b"
\ar "c";"b"
\ar "e";"d"
\ar "f";"d"
}$$
\end{exam}

The $m$-th associahedron is defined to be the polytope which has one vertex for every fully parenthesized 
word of length $m$.  Two vertices $(W,P_i)$ and $(W,P_j)$ are on the same
$k$-dimensional face if they share at least $m-k$ parentheses, i.e. $P_i\cap P_j$ has cardinality at least $m-k$.  
Thus our poset $\scK_m$ is exactly the face poset of the $m$-th associahedron, and so the geometric realization 
of the nerve of $\scK_m$ is simply the barycentric subdivision of the $m$-th associahedron.
Note that we are using a Fraktur font to distinguish the poset $\scK_m$ from the associahedron $K_m$ which is the 
geometric realization of its nerve (as a topological space). 

 The following lemma will prove to be surprisingly useful:
 
\begin{lem}\label{surprise}
Let $(W_k,P)<(W_k,P')$ in  $\scK_k$. Then the subposet
$$[(W_k,P),(W_k,P')]=\{(W_k,P'')\in\scK_k | (W_k,P)\le(W_k,P'')\le(W_k,P')\}$$
is isomorphic to the poset $\scI^m$ where $\scI$ is the poset $1< 0$ and $m$ is the number of parentheses in 
$(W_k,P)$ which are not in $(W_k,P')$.
In other words, the factorizations of a fixed morphism in $\scK_k$ form a commutative cubical diagram.
\end{lem}

\begin{proof}
We can uniquely associate to each element $(W_k,P'')$ in $[(W_k,P),(W_k,P')]$ a characteristic function on the 
set of parentheses in $(W_k,P)$ which are not in $(W_k,P')$ by giving the value 1 to each
parenthesis which occurs in $(W_k,P'')$ and 0 to any which do not so occur. But such a characteristic function 
is evidently the same thing as an object of $\scI^m$ and it is clear that
order relations match.
\end{proof}

If we take $(W_k,P')=(W_k,\emptyset)$, then geometrically this gives a decomposition of the associahedra into cubes.  
The decomposition of $K_4$ into 5 squares looks like this:

$$
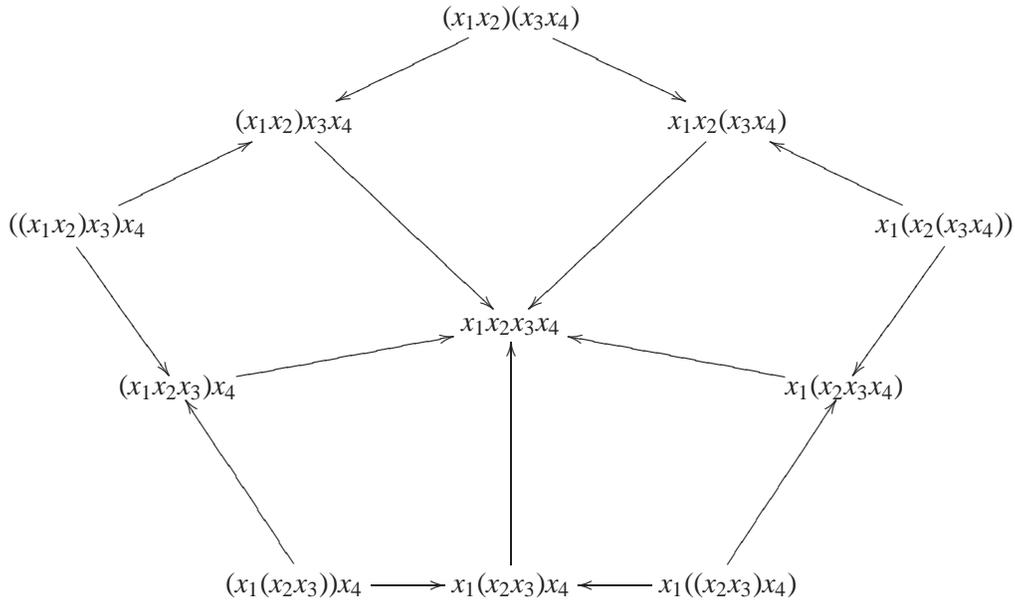
\xymatrix{
&&(x_1x_2)(x_3x_4)\ar[dl]\ar[dr]\\
&(x_1x_2)x_3x_4\ar[ddr] &&x_1x_2(x_3x_4)\ar[ddl]\\
{((x_1x_2)x_3)x_4}\save[]+<0pt,-8pt>*{}="a"\restore\ar[ur] && &&x_1(x_2(x_3x_4))\save[]+<0pt,-8pt>*{}="f"\restore\ar[ul]\\
&&x_1x_2x_3x_4\\
{}\save[]+<38pt,10pt>*{(x_1x_2x_3)x_4}="b"\ar[urr]\restore && &&{}\save[]+<-38pt,10pt>*{x_1(x_2x_3x_4)}="d"\ar[ull]\restore\\
\\
&\save[]+<0pt,8pt>*{}="c"\restore(x_1(x_2x_3))x_4\ar[r] &x_1(x_2x_3)x_4\ar[uuu] 
&\save[]+<0pt,8pt>*{}="e"\restore x_1((x_2x_3)x_4)\ar[l]
\ar "a";"b"
\ar "c";"b"
\ar "e";"d"
\ar "f";"d"
}$$

The decomposition of $K_5$ into 14 cubes can be found here:\newline
\url{http://arxiv.org/src/1005.3979v4/anc/cubical.flv}.
 
There is a one-to-one correspondence between parenthesized words and \textit{stable rooted trees}. Briefly
 these are planar rooted trees where each node has at least two input edges.
We refer to \cite{Lein} for a formal definition. The correspondence is given by labelling the leaves of such a
 tree with the labels $x_1, x_2, \dots x_n$ in left to right order. [In the
degenerate cases $n=1$ and $n=0$, the identity $id\in \scK_1$ corresponds to the tree with a single edge and no nodes and
 $\emptyset\in \scK_0$ corresponds to the empty tree with no edges and no nodes.]
Then for each node of the tree, except for the bottom root node, one takes the set of labels sitting over that
 node as one of the intervals $p_i\in P$ in the collection $P$, thus giving us
a parenthesized word $(W_n,P)$.  For example, here are all of the parenthesized words on the linear order $W_4$ 
and their corresponding stable rooted trees:

\begin{exam}
The following trees\newline
$\xymatrix@=5pt@M=-1pt@W=-1pt{
\ar@{-}[dr] &&\ar@{-}[dl]&&\ar@{-}[dr] &&\ar@{-}[dl]\\
&{\scriptscriptstyle\bullet}\ar@{-}[ddrr]&&&&{\scriptscriptstyle\bullet}\ar@{-}[ddll]\\
\\
&&&{\scriptscriptstyle\bullet}\ar@{-}[d]\\
&&&&&{}
}$\raisebox{-20pt}{,}
$\xymatrix@=5pt@M=-1pt@W=-1pt{
\ar@{-}[dr] &&\ar@{-}[dl]&&&\ar@{-}[ddll] &&&\ar@{-}[dddlll]\\
&{\scriptscriptstyle\bullet}\ar@{-}[drr]\\
&&&{\scriptscriptstyle\bullet}\ar@{-}[drr]\\
&&&&&{\scriptscriptstyle\bullet}\ar@{-}[d]\\
&&&&&
}$\raisebox{-20pt}{,}
$\xymatrix@=5pt@M=-1pt@W=-1pt{
\ar@{-}[ddrrr] &&&\ar@{-}[dr]&&\ar@{-}[dl] &&&\ar@{-}[dddllll]\\
&&&&{\scriptscriptstyle\bullet}\ar@{-}[dl]\\
&&&{\scriptscriptstyle\bullet}\ar@{-}[dr]\\
&&&&{\scriptscriptstyle\bullet}\ar@{-}[d]\\
&&&&
}$\raisebox{-20pt}{,}
$\xymatrix@=5pt@M=-1pt@W=-1pt{
\ar@{-}[dddrrrr] &&&\ar@{-}[dr]&&\ar@{-}[dl] &&&\ar@{-}[ddlll]\\
&&&&{\scriptscriptstyle\bullet}\ar@{-}[dr]\\
&&&&&{\scriptscriptstyle\bullet}\ar@{-}[dl]\\
&&&&{\scriptscriptstyle\bullet}\ar@{-}[d]\\
&&&&
}$\raisebox{-20pt}{,}
$\xymatrix@=5pt@M=-1pt@W=-1pt{
\ar@{-}[dddrrrr] &&&\ar@{-}[ddrr]&&\ar@{-}[dr] &&&\ar@{-}[dll]\\
&&&&&&{\scriptscriptstyle\bullet}\ar@{-}[dl]\\
&&&&&{\scriptscriptstyle\bullet}\ar@{-}[dl]\\
&&&&{\scriptscriptstyle\bullet}\ar@{-}[d]\\
&&&&
}$\raisebox{-20pt}{,}
$\xymatrix@=5pt@M=-1pt@W=-1pt{
\ar@{-}[ddr]&\ar@{-}[dd]&\ar@{-}[ddl] &&&\ar@{-}[dddlll]\\
\\
&{\scriptscriptstyle\bullet}\ar@{-}[dr]\\
&&{\scriptscriptstyle\bullet}\ar@{-}[d]\\
&&
}$\raisebox{-20pt}{,}
$\xymatrix@=5pt@M=-1pt@W=-1pt{
\ar@{-}[dddrrr]&&&\ar@{-}[ddr]&\ar@{-}[dd]&\ar@{-}[ddl]\\
\\
&&&&{\scriptscriptstyle\bullet}\ar@{-}[dl]\\
&&&{\scriptscriptstyle\bullet}\ar@{-}[d]\\
&&&
}$\raisebox{-20pt}{,}
$\xymatrix@=5pt@M=-1pt@W=-1pt{
\ar@{-}[dr]&&\ar@{-}[dl]&\ar@{-}[dddl]&&&\ar@{-}[dddllll]\\
&{\scriptscriptstyle\bullet}\ar@{-}[ddr]\\
\\
&&{\scriptscriptstyle\bullet}\ar@{-}[d]\\
&&
}$\raisebox{-20pt}{,}
$\xymatrix@=5pt@M=-1pt@W=-1pt{
\ar@{-}[dddrrr]&&\ar@{-}[dr]&&\ar@{-}[dl]&&\ar@{-}[dddlll]\\
&&&{\scriptscriptstyle\bullet}\ar@{-}[dd]\\
\\
&&&{\scriptscriptstyle\bullet}\ar@{-}[d]\\
&&&
}$\raisebox{-20pt}{,}
$\xymatrix@=5pt@M=-1pt@W=-1pt{
\ar@{-}[dddrrrr]&&&\ar@{-}[dddr]&&\ar@{-}[dr]&&\ar@{-}[dl]\\
&&&&&&{\scriptscriptstyle\bullet}\ar@{-}[ddll]\\
\\
&&&&{\scriptscriptstyle\bullet}\ar@{-}[d]\\
&&&&
}$\raisebox{-20pt}{,}
$\xymatrix@=5pt@M=-1pt@W=-1pt{
\ar@{-}[dddrrr]&&\ar@{-}[dddr]&&\ar@{-}[dddl]&&\ar@{-}[dddlll]\\
\\
\\
&&&{\scriptscriptstyle\bullet}\ar@{-}[d]\\
&&&
}$
represent $(x_1x_2)(x_3x_4)$, $((x_1x_2)x_3)x_4$, $(x_1(x_2x_3))x_4$, $x_1((x_2x_3)x_4)$, $x_1(x_2(x_3x_4))$, 
$(x_1x_2x_3)x_4$, $x_1(x_2x_3x_4)$, $(x_1x_2)x_3x_4$,
$x_1(x_2x_3)x_4$, $x_1x_2(x_3x_4)$, $x_1x_2x_3x_4$ respectively.

\end{exam}

The poset structure on $\scK_n$ of Definition \ref{Kmdefi} can be described in terms of trees as follows:
 $T<T'$ if $T'$ can be obtained from $T$ by shrinking some of the edges of $T$.
It is of course more convenient to use parenthesized words when describing $\scK_n$ as a poset. However the 
language of trees is more convenient to describe the
operad structure on the $\scK_n$.

The $\{\scK_i\}_{i\ge 0}$ form an operad in $CAT$, the category of small categories.  Given stable rooted trees
 $S \in \scK_m$ and $T_i \in \scK_{k_i}$ for $i = 1,2, \dots, m$,
we obtain a new tree by grafting the root of $T_i$ to the $i$-th leaf of $S$.  In terms of parenthesized 
words we are  substituting the word for $T_i$ in place of the $i$-th character of the word for $S$, and 
reindexing to insure that all characters in the resulting word are distinct.  This only makes sense if $m\ge 2$ 
and all $k_i\ge 2$.  If $S=id\in \scK_1$, we define the composed tree to be  $T_1$. If $T_i=id\in \scK_1$, then we
 leave the $i$-th leaf of $S$ unchanged. If $T_i=0\in \scK_0$, then we delete the $i$-th leaf of $S$. 
If this leaves only one input edge for the node below, we delete that node as well.
If it leaves no input edges for the node below, we delete both that node and the edge below.  We apply this
 algorithm recursively: if the next node below receives only one input edge or no
input edges, we delete that node or that node together with the edge below, and so on.  In the special case
 when $\sum k_i = 1$ or $\sum k_i=0$,
the resulting degenerate trees are defined to be $id$ or $0$ respectively.

This process is is clearly functorial in each of $\scK_m, \scK_{k_1}, \scK_{k_2}, \dots \scK_{k_m}$, and 
so we obtain a functor 
$$\gamma_{m,k_1, k_2, \dots k_m}: \scK_m \times \prod_1^m \scK_{k_i} \to \scK_{\sum_1^m k_i}.$$
These functors define a categorical operad $\scK = \{\scK_i\}_{i\ge 0}$.

The associahedral operad $\scK=\{\scK_i\}_{i\ge 0}$ has an operadic filtration 
$$\scK^{(2)}\subset \scK^{(3)}\subset \scK^{(4)}\dots,$$
where $\scK^{(n)}_i$ is the subposet of $\scK_i$
consisting of trees where each node has input valence $\le n$ (i.e has at most $n$ incoming edges). We 
define $\scK^{(\infty)}=\scK$. We note for future reference that if
an element $(W,P)$ of $\scK_k$ lies in filtration $n$ and $(W',P')<(W,P)$ then $(W',P')$ also lies in 
filtration $n$. This follows from our
description above of the poset structure in terms of trees. 

\begin{propo}\label{connectivity} The poset $\scK^{(n)}_i$ is the face poset of a subcomplex of the (unsubdivided) 
associahedron $K_i$. This subcomplex contains all cells of
$K_i$ of dimension $\le n-2$. Consequently the nerve of $\scK^{(n)}_i$ is $(n-3)$-connected. In particular
 if $n\ge 4$ the nerve of $\scK^{(n)}_i$
is simply connected.
\end{propo}

\begin{proof} We note that if an element $(W,P)$ of $\scK_i$ is in $\scK^{(n)}_i$ and $(W',P')<(W,P)$ 
then $(W',P')$ is also contained in
$\scK^{(n)}_i$, since the tree representing  $(W,P)$ is obtained from the tree representing $(W',P')$ by 
shrinking internal edges. It follows
that $\scK^{(n)}_i$ is the face poset of a subcomplex of $K_i$. Now the vertices of $K_i$ are parametrized
 by the elements of $\scK^{(2)}_i$,
which are represented by binary trees. It follows that the cells of $K_i$ of dimension $j$ are obtained by
 shrinking $j$ internal edges of a binary
tree. It easily follows that the dimension of the cell parametrized by a given tree is the sum over all
nodes of the incoming valence of that node minus 2. Thus the maximal possible incoming valence of a node in 
a tree parametrizing a cell of dimension $j$ is $j+2$.
Hence the subcomplex of $K_i$ parametrized by $\scK^{(n)}_i$ contains all cells of
$K_i$ of dimension $\le n-2$.

Now the nerve of $\scK^{(n)}_i$ is the barycentric subdivision of this subcomplex of $K_i$. Moreover $K_i$
 is obtained from this subcomplex
by adding cells of dimensions $\ge n-1$. Since $K_i$ is contractible, it follows that the complex and hence
 the nerve of $\scK^{(n)}_i$
is $(n-3)$-connected.
\end{proof}

\begin{rema} $\scK^{(n)}_i$ is generally larger than the face poset of the $(n-2)$-skeleton of $K_i$. For 
instance $\scK^{(3)}_5$ is the
face poset of the subcomplex of the 3-dimensional associahedron $K_5$ consisting of all the edges together 
with the three square faces.
\end{rema}

\begin{rema}
Our categorical operad $\scK$ is almost the same as Leinster's $StTr$ (\cite[pages 233-234]{Lein}). The only difference is that
he has $StTr(0)=\emptyset$, whereas we have $\scK_0=\{0\}$. So our approach encodes the notion of a unit for algebras over $\scK$.
Leinster expected that the nerve of $StTr(k)=\scK_k$ is homeomorphic to the 
associahedron , which we prove.  Thus Leinster's topological operad is precisely the same as Stasheff's.\\
The tree description of a $CAT$-operad containing $\scK$ appears in \cite{GK}.
\end{rema}

\begin{rema}
Since the nerve of a product in $CAT$ is a product in $TOP$, it follows that the nerve of a $\scK^{(n)}$ 
algebra is an $A_n$-space
in the sense of Stasheff.
\end{rema}

\section{$A_n$-monoidal categories and coherence}

\begin{defi}
For $n=2,3,\dots,\infty$, an $A_n$-monoidal category is a category $\scC$ together with multiplications
 $\mu_k: \scC^k \rightarrow \scC$ for $0\le k<n+1$ such that
\begin{enumerate}
\item $\mu_1:  \scC \rightarrow \scC$ is the identity functor.
\item $\mu_0: * \rightarrow \scC$ is an object $0\in\scC$ that acts as a strict unit in the sense that
$$\mu_k(Id_\scC^{\,i}\times 0\times Id_\scC^{\,j}) = \mu_{k-1}$$
for any $i,j$ such that $i+j=k-1$.
\end{enumerate}
$\scC$ is also equipped with natural transformations (associators) 
$$\alpha^{i,j,k}: \mu_{i+1+k}\circ(Id_\scC^{\,i} \times \mu_j \times Id_\scC^{\,k} ) \longrightarrow 
\mu_{i+j+k},$$
for $0\le i+j+k<n+1$, satisfying
\begin{eqnarray*}
\lefteqn{(i)}\qquad\qquad
&\alpha^{i,0,k}, \alpha^{i,1,k} \mbox{ and } \alpha^{0,j,0} \mbox{ are the identity}\qquad\qquad\qquad
\qquad\qquad\qquad\qquad\\
\end{eqnarray*}
and the coherence conditions specified by the following commutative diagrams
\begin{eqnarray*}
\lefteqn{\raisebox{-20pt}{\mbox{(ii)}}}
&\qquad&\qquad\qquad\xymatrix{
\mu_{a+b+d+2}(\overline{A},0,\overline{B},\mu_c(\overline{C}),\overline{D})
\ar@{=}[r]\ar[d]^{\alpha^{a+b+1,c,d}_{(\overline{A},0,\overline{B}),\overline{C},\overline{D}}}
&\mu_{a+b+d+1}(\overline{A},\overline{B},\mu_c(\overline{C}),\overline{D})
\ar[d]^{\alpha^{a+b,c,d}_{(\overline{A},\overline{B}),\overline{C},\overline{D}}}\\
\mu_{a+b+c+d+1}(\overline{A},0,\overline{B},\overline{C},\overline{D})\ar@{=}[r]
&\mu_{a+b+c+d}(\overline{A},\overline{B},\overline{C},\overline{D})
}\\
\lefteqn{\raisebox{-20pt}{\mbox{(iii)}}}
&\qquad&\qquad\qquad\xymatrix{
\mu_{a+c+d+2}(\overline{A},\mu_b(\overline{B}),\overline{C},0,\overline{D})
\ar@{=}[r]\ar[d]^{\alpha^{a,b,c+d+1}_{\overline{A},\overline{B},(\overline{C},0,\overline{D})}}
&\mu_{a+c+d+1}(\overline{A},\mu_b(\overline{B}),\overline{C},\overline{D})
\ar[d]^{\alpha^{a,b,c+d}_{\overline{A},\overline{B},(\overline{C},\overline{D})}}\\
\mu_{a+b+c+d+1}(\overline{A},\overline{B},\overline{C},0,\overline{D})\ar@{=}[r]
&\mu_{a+b+c+d}(\overline{A},\overline{B},\overline{C},\overline{D})
}\\
\lefteqn{\raisebox{-20pt}{\mbox{(iv)}}}
&\qquad&\qquad\qquad\xymatrix{
\mu_{a+d+1}(\overline{A},\mu_{b+c+1}(\overline{B},0,\overline{C}),\overline{D})
\ar@{=}[r]\ar[d]^{\alpha^{a,b+c+1,d}_{\overline{A},(\overline{B},0,\overline{C}),\overline{D}}}
&\mu_{a+d+1}(\overline{A},\mu_{b+c}(\overline{B},\overline{C}),\overline{D})
\ar[d]^{\alpha^{a,b+c,d}_{\overline{A},(\overline{B},\overline{C}),\overline{D}}}\\
\mu_{a+b+c+d+1}(\overline{A},\overline{B},0,\overline{C},\overline{D})\ar@{=}[r]
&\mu_{a+b+c+d}(\overline{A},\overline{B},\overline{C},\overline{D})
}\\
\end{eqnarray*}
\begin{eqnarray*}
\lefteqn{\raisebox{-20pt}{\mbox{(v)}}}
&\qquad&\xymatrix{
\mu_{a+c+e+2}(\overline{A},\mu_b(\overline{B}),\overline{C},\mu_d(\overline{D}),\overline{E})
\ar[d]^{\alpha^{a+c+1,d,e}_{(\overline{A},\mu_b(\overline{B}),\overline{C}),\overline{D},\overline{E}}}
\ar[rr]^{\alpha^{a,b,c+e+1}_{\overline{A},\overline{B},(\overline{C},\mu_d(\overline{D}),\overline{E})}}
&&\mu_{a+b+c+e+1}(\overline{A},\overline{B},\overline{C},\mu_d(\overline{D}),\overline{E})
\ar[d]^{\alpha^{a+c+1,d,e}_{(\overline{A},\overline{B},\overline{C}),\overline{D},\overline{E}}}\\
\mu_{a+c+d+e+1}(\overline{A},\mu_b(\overline{B}),\overline{C},\overline{D},\overline{E})
\ar[rr]^{\alpha^{a,b,c+d+e}_{\overline{A},\overline{B},(\overline{C},\overline{D},\overline{E})}}
&&\mu_{a+b+c+d+e}(\overline{A},\overline{B},\overline{C},\overline{D},\overline{E})
}\\
\lefteqn{\raisebox{-20pt}{\mbox{(vi)}}}
&\qquad&\xymatrix{
\mu_{a+e+1}(\overline{A},\mu_{b+d+1}(\overline{B},\mu_c(\overline{C}),\overline{D}),\overline{E})
\ar[rrr]^{\mu_{a+e+1}(id_{\overline{A}},\alpha^{b,c,d}_{\overline{B},\overline{C},\overline{D}},id_{\overline{E}})}
\ar[d]^{\alpha^{a,b+d+1,e}_{\overline{A},(\overline{B},\mu_c(\overline{C}),\overline{D}),\overline{E}}}
&&&\mu_{a+e+1}(\overline{A},\mu_{b+c+d}(\overline{B},\overline{C},\overline{D}),\overline{E})
\ar[d]^{\alpha^{a,b+c+d,e}_{\overline{A},(\overline{B},\overline{C},\overline{D}),\overline{E}}}\\
\mu_{a+b+d+e+1}(\overline{A},\overline{B},\mu_c(\overline{C}),\overline{D},\overline{E})
\ar[rrr]^{\alpha^{a+b,c,d+e}_{(\overline{A},\overline{B}),\overline{C},(\overline{D},\overline{E})}}
&&&\mu_{a+b+c+d+e}(\overline{A},\overline{B},\overline{C},\overline{D},\overline{E})
}
\end{eqnarray*}
\end{defi}
Here $\overline{A}$, $\overline{B}$, $\overline{C}$, $\overline{D}$, $\overline{E}$ are taken to be objects
of $\scC^a$, $\scC^b$, $\scC^c$, $\scC^d$, $\scC^e$, respectively.

Essentially coherence conditions (i)-(iv) require the associators to be compatible with the strict unit 0, 
while (v) and (vi) just say that if we are removing two pairs of matching 
parentheses in a multiplication,
it doesn't matter which we remove first.

\begin{rema}
In \cite[pages 93-94]{Lein}, Leinster defines the notion of a lax monoidal category, which is similar in spirit to the above
definition, but there are some crucial differences. A lax monoidal category in his sense, has multiplications
$$\mu_k(A_1,A_2,\dots,A_k)=(A_1\otimes A_2\otimes\dots\otimes A_k)$$
for all $k\in \mathbb{N}$ together with natural transformations
$$
\gamma^{k_1,\ldots,k_n}: \mu_n\circ (\mu_{k_1}\otimes\ldots\otimes \mu_{k_n})\to \mu_{k_1+\ldots +k_n}
$$
and a natural transformation
$$\iota_A: A\longrightarrow \mu_1(A)=(A).$$
The natural transformations $\gamma$ satisfy a coherence condition which is essentially our coherence conditions (v) and (vi) combined
into a single diagram. There is no unit condition for $\mu_0$ (so one might as well require the existence of $\mu_k$ for
$k>0$ only). Moreover his natural transformation
$\iota$ is not the identity. Thus a lax monoidal category in his sense possesses arbitrarily long nondegenerate strings of composable natural
transformations between unary multiplications
$$A\stackrel{\iota_{A}}{\longrightarrow}(A)\stackrel{\iota_{(A)}}{\longrightarrow}((A))\stackrel{\iota_{((A))}}{\longrightarrow}\dots$$
It follows that the operad controlling such a structure has an infinite dimensional nerve.
\end{rema}

The main result of this paper is:

\begin{theo}\label{main}
A category $\scC$ is a $\scK^{(n)}$-algebra iff it is an $A_n$-monoidal category.
\end{theo}

\begin{proof}
Given an action, $\theta_i: \scK^{(n)}_i\times\scC^i\longrightarrow\scC$, define $\mu_i:\scC^i\longrightarrow\scC$ 
to be the restriction
of this action to $\{x_1x_2\dots x_i\}\times\scC^i$, where $x_1x_2\dots x_i=(W_i,\emptyset)$ is the terminal object
 of $\scK_i$. This makes sense for
$0\le i<n+1$, since in those cases $x_1x_2\dots x_i$ is contained in the $n$-th filtration $\scK^{(n)}$. We 
then define $\alpha^{i,j,k}$ to be
the restriction of $\theta_{i+j+k}$ to $\left\{\left(W_{i+j+k},\{[x_{i+1},x_{i+j}]\}\right)
\longrightarrow(W_{i+j+k},\emptyset)\right\}\times\scC^{i+j+k}$,
for $0\le i+j+k<n+1$.

Conditions (1), (2) and (i) follow from the fact that $(x_1)\in\scK^{(n)}_1$ is the identity of the operad 
and composing the constant $0\in\scK^{(n)}_0$
into any input of $\{x_1x_2\dots x_i\}\in\scK^{(n)}_i$ gives $\{x_1x_2\dots x_{i-1}\}\in\scK^{(n)}_{i-1}$. 
Conditions (ii)-(iv) also follow from the latter fact.
Finally conditions (v) and (vi) follow from the restriction of $\theta_{a+b+c+d+e}$ to 
$\scD\times\scC^{a+b+c+d+e}$ and $\scD'\times\scC^{a+b+c+d+e}$,
where $\scD$ and $\scD'$ are the following commutative diagrams in $\scK^{(n)}_{a+b+c+d+e}$:
$$\xymatrix{
\left(W_{a+b+c+d+e},\left\{[x_{a+1},x_{a+b}],[x_{a+b+c+1},x_{a+b+c+d}]\right\}\right)\ar[r]\ar[d]
&\left(W_{a+b+c+d+e},\left\{[x_{a+1},x_{a+b}]\right\}\right)\ar[d]\\
\left(W_{a+b+c+d+e},\left\{[x_{a+b+c+1},x_{a+b+c+d}]\right\}\right)\ar[r]
&(W_{a+b+c+d+e},\emptyset)
}$$
$$\xymatrix{
\left(W_{a+b+c+d+e},\left\{[x_{a+1},x_{a+b+c+d}],[x_{a+b+1},x_{a+b+c}]\right\}\right)\ar[r]\ar[d]
&\left(W_{a+b+c+d+e},\left\{[x_{a+b+1},x_{a+b+c}]\right\}\right)\ar[d]\\
\left(W_{a+b+c+d+e},\left\{[x_{a+1},x_{a+b+c+d}]\right\}\right)\ar[r]
&(W_{a+b+c+d+e},\emptyset)
}$$
respectively.

Conversely suppose that $\scC$ is an $A_n$-monoidal category. Then we define
$$\theta_i:\mbox{Obj}(\scK^{(n)}_i)\times\scC^i\longrightarrow\scC$$
by induction on $i$ as follows.  We define $\theta_0$ to be $\mu_0$ and $\theta_1$ to be $\mu_1=id_\scC$. 
Having defined $\theta_j$ for $j<i$,
consider an object $T$ in $\scK^{(n)}_i$ represented by a tree
$$\xymatrix@=5pt@M=-1pt@W=-1pt{
T_{1} &T_{2} &\dots &\ T_{k}\\
\ar@{-}[ddrr]&\ar@{-}[ddr]& &\ar@{-}[ddl]\\
\\
&&{\scriptscriptstyle\bullet}\ar@{-}[dd]\\
\\
&&}$$
with $k<n+1$ and where $T_j$ has $m_j$ input edges, so that $m_1+m_2+\dots+m_k=i$. Let 
$(\overline{A}_1,\overline{A}_2,\dots,\overline{A}_k)$
represent an object in $\scC^i$, with $\overline{A}_j\in\mbox{Obj}\left(\scC^{m_j}\right)$, $j=1,2,\dots,k$. By induction
$\theta_{m_j}(T_j,\overline{A}_j)$ are already defined for $j=1,2,\dots,k$.  We then define
$$\theta_i\left(T,\overline{A}_1,\overline{A}_2,\dots,\overline{A}_k)=
\mu_k(\theta_{m_1}(T_1,\overline{A}_1),\theta_{m_2}(T_2,\overline{A}_2),\dots,\theta_{m_k}(T_k,\overline{A}_k)\right)$$
[Here we use implicitly the canonical isomorphisms between the associahedral posets based on subintervals of $W_i$ with
the associahedral posets based on the standard linear orders $W_{m_i}$ of the same cardinality, c.f. \ref{Kmdefi}.]
We define $\theta_i$ for morphisms in $\scC^i$ similarly. This completes the induction.

Next we extend the definition of $\theta_i$ to define natural transformations
$$\theta_i:\mbox{IMor}(\scK^{(n)}_i)\times\mbox{Obj}(\scC^i)\longrightarrow\mbox{Mor}(\scC)$$
where $\mbox{IMor}(\scK^{(n)}_i)$ are the indecomposable morphisms in $\scK^{(n)}_i$, i.e. morphisms
 which can't be factored nontrivially
(or equivalently morphisms given by dropping a single pair of matching parentheses in a parenthesized word).
Again we proceed by induction on $i$, starting with $i=0$ and $i=1$ where these are vacuously defined. 
Now consider an indecomposable
morphism $\lambda:T\to T'$ in $\scK^{(n)}_i$, where $T$ has the form
$$
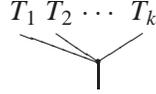
\xymatrix@=5pt@M=-1pt@W=-1pt{
T_{1} &T_{2} &\dots &\ T_{k}\\
\ar@{-}[ddrr]&\ar@{-}[ddr]& &\ar@{-}[ddl]\\
\\
&&{\scriptscriptstyle\bullet}\ar@{-}[dd]\\
\\
&&}$$
with $k<n+1$ and where $T_j$ has $m_j$ input edges, so that $m_1+m_2+\dots+m_k=i$. Then $\lambda$ is 
obtained by shrinking a single interior edge
in $T$. There are two possibilities: (1) an interior edge of some tree $T_j$ is shrunk or (2) an edge
 below some $T_j$ is shrunk. Now let
$(\overline{A}_1,\overline{A}_2,\dots,\overline{A}_k)$ represent an object in $\scC^i$, with
$\overline{A}_j\in\mbox{Obj}\left(\scC^{m_j}\right)$, $j=1,2,\dots,k$. In the first case we define

$\theta_i(\lambda,\overline{A}_1,\overline{A}_2,\dots,\overline{A}_k)=$\\
$\mu_k\left(id_{\theta_{m_1}(T_1,\overline{A}_1)},id_{\theta_{m_2}(T_2,\overline{A}_2)},\dots,
id_{\theta_{m_{j-1}}(T_{j-1},\overline{A}_{j-1})},\theta_{m_j}(\lambda',\overline{A_j}),
id_{\theta_{m_{j+1}}(T_{j+1},\overline{A}_{j+1})},\dots,id_{\theta_{m_k}(T_k,\overline{A}_k)}\right)
$

where $\lambda'$ is the indecomposable morphism in $\scK^{(n)}_{m_j}$ given by shrinking that particular edge. 
 In the second case we define
$$\theta_i(\lambda,\overline{A}_1,\overline{A}_2,\dots,\overline{A}_k)
=\alpha^{m_1+m2+\dots m_{j-1},m_j,m_{j+1}+m_{j+2}+\dots+m_k}_{(\overline{A}_1,\overline{A}_2,\dots,
\overline{A}_{j-1}),\overline{A}_j,(\overline{A}_{j+1},\dots,\overline{A}_k)}.$$
This completes the induction.

Finally to extend $\theta_i$ to all morphisms in $\scK_i$, we must show that for any factorization of a 
morphism in $\scK^{(n)}_i$ into indecomposable morphisms,
the corresponding composition of natural transformations defines the same morphism in $\scC$. But according 
to Lemma \ref{surprise}, the factorizations of any morphism
in $\scK_i$ give rise to a cubical diagram in $\scC$. According to coherence conditions (v) and (vi) of an 
$A_n$-monoidal category, all the 2-dimensional faces of this cubical
diagram commute. It is an elementary consequence that the entire cubical diagram in $\scC$ commutes, c.f. 
Lemma \ref{CubicalLemma} below.  It follows that there are well defined functors:
$$\theta_i: \scK^{(n)}_i\times\scC^i\longrightarrow\scC$$
for all $i\ge 0$.  The fact that $\theta_i$ are compatible with the operadic compositions
$$\gamma_{m,k_1, k_2, \dots k_m}: \scK^{(n)}_m \times \prod_1^m \scK^{(n)}_{k_i} \to \scK^{(n)}_{\sum_1^m k_i}$$
follows from the inductive construction of $\theta_i$ if all the $k_i>1$. If $k_i\le 1$ or $m=1$, the 
compatibility follows from conditions (1), (2) and (i)-(iv) of the
definition of an $A_n$-monoidal category.
\end{proof}

\begin{lem}\label{CubicalLemma} A cubical diagram in any category commutes iff each of its 2-dimensional faces commutes.
\end{lem}

\begin{proof} We proceed by induction on the dimension of the cube. The statement is vacuously true if the dimension is $\le 2$.
Suppose it is true for all cubical diagrams of dimension $<m$, and suppose we are given an $m$-dimensional cubical diagram.
Consider two edge paths from the initial object $A$ of the diagram to $Z$, the terminal object. Let these edge paths factor as
$$A\stackrel{\alpha}{\longrightarrow}B\stackrel{f}{\longrightarrow}Z,\qquad 
A\stackrel{\beta}{\longrightarrow}C\stackrel{g}{\longrightarrow}Z$$
respectively, where $\alpha$ and $\beta$ are edges of the diagram and $f$ and $g$ are composites of the remainders of these
edge paths. If $\alpha=\beta$, then by induction $f=g$ and we are done. Otherwise let
$$\xymatrix{A\ar[r]^\alpha\ar[d]^\beta &B\ar[d]^\gamma\\
C\ar[r]^\delta &D
}$$
be the 2-dimensional face spanned by $\alpha$ and $\beta$. Pick any edge path $h:D\longrightarrow Z$, and consider the
diagram
$$
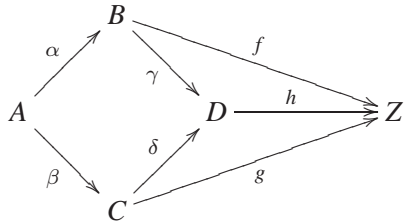
\xymatrix{
&B\ar[dr]_\gamma\ar[drrr]^f\\
A\ar[ur]^\alpha\ar[dr]_\beta &&D\ar[rr]^{h\quad} &&Z\\
&C\ar[ur]^\delta\ar[urrr]_g
}$$
By hypothesis the square commutes and by induction the two triangles commute. Hence $f\alpha=g\beta$. This completes the
induction and proof.
\end{proof}

\section{Relation to coherence theorems for monoidal categories}

\begin{defi} We say that an $A_n$-monoidal category is {\it undirected\/} if all the associativity natural
 transformations $\alpha^{i,j,k}$ are isomorphisms.
\end{defi}

\begin{propo} An undirected $A_n$-monoidal category is a monoidal category if $n\ge 4$.
\end{propo}

\begin{proof} If $\scC$ is an undirected $A_n$-monoidal category, then the corresponding action functors 
$\theta_i:\scK^{(n)}_i\times\scC^i\longrightarrow\scC$
extend to $\theta_i:\overline{\scK}^{(n)}_i\times\scC^i\longrightarrow\scC$, where $\overline{\scK}^{(n)}_i$ 
is obtained from $\scK^{(n)}_i$ by formally
inverting all the morphisms. By Proposition \ref{connectivity} the nerve of $\scK^{(n)}_i$ is simply connected.  
Now recalling that inverting all the morphisms
in a connected category has the effect of killing off the higher homotopy groups of its nerve 
(c.f. \cite[Proposition 1]{Q}), we see that
the nerve of $\overline{\scK}^{(n)}_i$
is contractible, and it follows that the objects of $\scK^{(2)}_i$ are connected to each other by 
uniquely defined isomorphisms in $\overline{\scK}^{(n)}_i$.
The images of these isomorphisms under $\theta_i$ specify uniquely defined natural isomorphisms connecting 
all possible different ways of associating the
binary product $\mu_2:\scC^2\longrightarrow\scC$ into an $i$-fold product $\scC^i\longrightarrow\scC$ 
so that all diagrams involving them commute. Thus
$\scC$ is a monoidal category, in the classical sense of Mac Lane.
\end{proof}

Next we derive LaPlaza's coherence theorem \cite{LaPl}, which generalizes Mac Lane's 
coherence theorem to the case
where the associativity natural transformation for a monoidal structure on a category is not 
required to be an isomorphism. We begin with a preliminary version of this result.

\begin{theo}\label{CohThm0}
Let $(\scC,\Box,0,\eta)$ be a directed monoidal category with a strict unit.  That is, $\Box:\scC\times\scC\longrightarrow\scC$
is a bifunctor and 0 is an object of $\scC$ which serves as a strict unit for $\Box$, i.e. the restrictions of $\Box$ to
$0\times\scC$ and $\scC\times 0$ are the identity. Finally $\eta_{A,B,C}:(A\Box B)\Box C\longrightarrow A\Box(B\Box C)$
is a natural transformation (not necessarily an isomorphism) such that $\eta_{A,B,C}$ is the identity whenever one of
$A$, $B$, $C$ is 0 and such that the pentagonal diagram
$$\xymatrix@=5pt@M=5pt@W=5pt{
&&(A\Box B)\Box(C\Box D)\ar[ddrr]^{\eta_{A,B,C\Box D}}\\
\\
((A\Box B)\Box C)\Box D\ar[uurr]^{\eta_{A\Box B,C,D}}\ar[ddddr]^{\eta_{A,B,C}\Box id_D} && &&A\Box(B\Box(C\Box D))\\
\\
\\
\\
&(A\Box(B\Box C))\Box D\ar[rr]^{\eta_{A,B\Box C,D}} &&A\Box((B\Box C)\Box D)\ar[uuuur]^{id_A\Box\eta_{B,C,D}}
}$$
commutes. Then $\scC$ can be endowed with the structure of an $A_\infty$-monoidal category.
\end{theo}

\begin{proof}
We define $\mu_0(*)=0$, $\mu_1$ to be the identity, $\mu_2=\Box$ and then we inductively define $\mu_i$ to be the
composite
$$\xymatrix{
\scC^i=\scC\times\scC^{i-1}\ar[rr]^{id_\scC\times\mu_{i-1}}&&\scC\times\scC\stackrel{\Box}{\longrightarrow}\scC
}$$
Thus $\mu_i(A_1,A_2,A_3,\dots,A_i)=A_1\Box(A_2\Box(A_3\Box(\dots\Box(A_{i-1}\Box A_i)\dots)))$ and we have
$$\xymatrix{
(*)
&&&\mu_{a+b}(\overline{A},\overline{B})=\mu_{a+1}(\overline{A},\mu_b(\overline{B}))&&&&&
}$$
for any objects $\overline{A}\in\scC^a$, $\overline{B}\in\scC^b$.

Now let $\overline{B}\in\scC^b$ and $\overline{C}\in\scC^c$. We define the associativity $\alpha^{0,b,c}_{\overline{B},\overline{C}}$
inductively on $b$. We assume $c>0$, since $\alpha^{0,b,0}$ is required by definition to be the identity. For $b=1$, we 
also require $\alpha^{0,1,c}$
 to be the identity. So suppose $b>1$ and $\overline{B}=(B_1,\overline{B}')$. Then we define $\alpha^{0,b,c}$ to
be the composite\newline
$\mu_{c+1}(\mu_b(\overline{B}),\overline{C})=\mu_2(\mu_b(\overline{B}),\mu_c(\overline{C}))=
(B_1\Box\mu_{b-1}(\overline{B}'))\Box\mu_c(\overline{C})$
$$\xymatrix{
&&\qquad\ar[rrrr]^{\eta_{B_1,\mu_{b-1}(\overline{B}'),\mu_c(\overline{C})}} &&&&B_1\Box(\mu_{b-1}(\overline{B}')\Box\mu_c(\overline{C}))
}$$
followed by the composite
$$B_1\Box(\mu_{b-1}(\overline{B}')\Box\mu_c(\overline{C}))=B_1\Box\mu_2(\mu_{b-1}(\overline{B}'),\mu_c(\overline{C}))
=B_1\Box\mu_{c+1}(\mu_{b-1}(\overline{B}'),\overline{C})$$
$$\xymatrix{
&&&\ar[rrr]^{id_{B_1}\Box\alpha^{0,b-1,c}_{\overline{B}',\overline{C}}\qquad\qquad}&&&B_1\Box\mu_{b+c-1}(\overline{B}',
\overline{C})=\mu_{b+c}(\overline{B},\overline{C}).
}$$
This completes the inductive definition of $\alpha^{0,b,c}_{\overline{B},\overline{C}}$.
We then define $\alpha^{a,b,c}_{\overline{A},\overline{B},\overline{C}}$ to be the composite\newline
$\mu_{a+c+1}(\overline{A},\mu_b(\overline{B}),\overline{C})=\mu_{a+1}(\overline{A},\mu_{c+1}(\mu_b(\overline{B}),\overline{C}))$
$$\xymatrix{
&&\ar[rrr]^{\mu_{a+1}(id_{\overline{A}},\alpha^{0,b,c}_{\mu_b(\overline{B}),\overline{C}})
\qquad\qquad\qquad\qquad}&&&\mu_{a+1}(\overline{A},\mu_{b+c}(\overline{B},\overline{C}))
=\mu_{a+b+c}(\overline{A},\overline{B},\overline{C}).
}$$
Note that this implies that $\alpha^{a,b,c}_{\overline{A},\overline{B},\overline{C}}$ is the identity if $c=0$, and that
$$\xymatrix{
(**)
&&&\alpha^{a_1+a_2,b,c}_{(\overline{A}_1,\overline{A}_2),\overline{B},\overline{C}}=
\mu_{a_1+1}(id_{A_1},\alpha^{a_2,b,c}_{\overline{A}_2,\overline{B},\overline{C}}).&&&&&
}$$

Conditions (1), (2), (i)-(iv) for an $A_\infty$-monoidal category are either true by construction or 
follow by a straight forward induction argument using the
hypotheses that 0 is a strict unit for $\Box$ and that $\eta_{A,B,C}$ is the identity whenever one of $A$, $B$ or $C$ is 0.

By (**) the verification of condition (v) reduces to the special case of the diagram
$$\xymatrix{
\mu_{c+e+2}(\mu_b(\overline{B}),\overline{C},\mu_d(\overline{D}),\overline{E})
\ar[d]^{\alpha^{c+1,d,e}_{(\mu_b(\overline{B}),\overline{C}),\overline{D},\overline{E}}}
\ar[rr]^{\alpha^{0,b,c+e+1}_{\overline{B},(\overline{C},\mu_d(\overline{D}),\overline{E})}}
&&\mu_{b+c+e+1}(\overline{B},\overline{C},\mu_d(\overline{D}),\overline{E})
\ar[d]^{\alpha^{c+1,d,e}_{(\overline{B},\overline{C}),\overline{D},\overline{E}}}\\
\mu_{c+d+e+1}(\mu_b(\overline{B}),\overline{C},\overline{D},\overline{E})
\ar[rr]^{\alpha^{0,b,c+d+e}_{\overline{B},(\overline{C},\overline{D},\overline{E})}}
&&\mu_{b+c+d+e}(\overline{B},\overline{C},\overline{D},\overline{E})
}$$
since the general diagram for (v) can be obtained from this one by applying the functor $\mu_{a+1}(\overline{A},-)$ to it.
By (*) and (**), this diagram in turn is the same as the diagram
$$\xymatrix{
\mu_{2}(\mu_b(\overline{B}),\mu_{c+e+1}(\overline{C},\mu_d(\overline{D}),\overline{E}))
\ar[d]^{\mu_2(id_{\mu_b(\overline{B})},\alpha^{c,d,e}_{\overline{C},\overline{D},\overline{E}})}
\ar[rrr]^{\alpha^{0,b,1}_{\overline{B},\mu_{c+e+1}(\overline{C},\mu_d(\overline{D}),\overline{E})}}
&&&\mu_{b+1}(\overline{B},\mu_{c+e+1}(\overline{C},\mu_d(\overline{D}),\overline{E}))
\ar[d]^{\mu_{b+1}(id_{\overline{B}},\alpha^{c,d,e}_{\overline{C},\overline{D},\overline{E}})}\\
\mu_{2}(\mu_b(\overline{B}),\mu_{c+d+e}(\overline{C},\overline{D},\overline{E}))
\ar[rrr]^{\alpha^{0,b,1}_{\overline{B},\mu_{c+d+e}(\overline{C},\overline{D},\overline{E})}}
&&&\mu_{b+1}(\overline{B},\mu_{c+d+e}(\overline{C},\overline{D},\overline{E}))
}$$
This last diagram in turn commutes because
$$\alpha^{0,b,1}_{\overline{B},X}:\mu_2(\mu_b(\overline{B}),X)\longrightarrow\mu_{b+1}(\overline{B},X)$$
is a natural transformation.

By similar reasoning, the verification of condition (vi) reduces to the special case of the diagram
$$\xymatrix{
\mu_{e+1}(\mu_{b+d+1}(\overline{B},\mu_c(\overline{C}),\overline{D}),\overline{E})
\ar[rrr]^{\mu_{e+1}(\alpha^{b,c,d}_{\overline{B},\overline{C},\overline{D}},id_{\overline{E}})}
\ar[d]^{\alpha^{0,b+d+1,e}_{(\overline{B},\mu_c(\overline{C}),\overline{D}),\overline{E}}}
&&&\mu_{e+1}(\mu_{b+c+d}(\overline{B},\overline{C},\overline{D}),\overline{E})
\ar[d]^{\alpha^{0,b+c+d,e}_{(\overline{B},\overline{C},\overline{D}),\overline{E}}}\\
\mu_{b+d+e+1}(\overline{B},\mu_c(\overline{C}),\overline{D},\overline{E})
\ar[rrr]^{\alpha^{b,c,d+e}_{\overline{B},\overline{C},(\overline{D},\overline{E})}}
&&&\mu_{b+c+d+e}(\overline{B},\overline{C},\overline{D},\overline{E})
}$$
By (*) it follows that this diagram is unchanged if we replace $E$ throughout by $\mu_e(E)$. Hence we
may as well suppose that $e=1$ and $\overline{E}=E$ is an object of $\scC$. Then by the inductive definition
of $\alpha^{0,i,j}$ and (**), we can factor the diagram as follows:
$$\xymatrix{
(B_1\Box\mu_{b+d}(\overline{B}',\mu_c(\overline{C}),\overline{D}))\Box E
\ar[rrr]^{(id_{B_1}\Box\alpha^{b-1,c,d}_{\overline{B}',\overline{C},\overline{D}})\Box id_{E}}
\ar[d]^{\eta_{B_1,\mu_{b+d}(\overline{B}',\mu_c(\overline{C}),\overline{D}),E}}
&&&(B_1\Box\mu_{b+c+d-1}(\overline{B}',\overline{C},\overline{D}))\Box E
\ar[d]^{\eta_{B_1,\mu_{b+c+d-1}(\overline{B}',\overline{C},\overline{D}),E}}\\
B_1\Box(\mu_{b+d}(\overline{B}',\mu_c(\overline{C}),\overline{D})\Box E)
\ar[rrr]^{id_{B_1}\Box(\alpha^{b-1,c,d}_{\overline{B}',\overline{C},\overline{D}}\Box id_{E})}
\ar[d]^{id_{B_1}\Box\alpha^{0,b+d,1}_{(\overline{B}',\mu_c(\overline{C}),\overline{D}),\overline{E}}}
&&&B_1\Box(\mu_{b+c+d-1}(\overline{B}',\overline{C},\overline{D})\Box E)
\ar[d]^{id_{B_1}\Box\alpha^{0,b+c+d-1,1}_{(\overline{B},\overline{C},\overline{D}),E}}
\\
B_1\Box\mu_{b+d+1}(\overline{B}',\mu_c(\overline{C}),\overline{D},E)
\ar[rrr]^{id_{B_1}\Box\alpha^{b-1,c,d+1}_{\overline{B},\overline{C},(\overline{D},E)}}
&&&B_1\Box\mu_{b+c+d}(\overline{B}',\overline{C},\overline{D},E)
}$$
The upper square commutes by naturality of $\eta$. The commutativity of the lower square corresponds
to a reduction of the problem from $b$ to $b-1$. Recursing on this reduction we reduce to the case $b=0$,
i.e. showing that the diagram
$$\xymatrix{
\mu_{d+1}(\mu_c(\overline{C}),\overline{D})\Box E
\ar[rrr]^{\alpha^{0,c,d}_{\overline{C},\overline{D}}\Box id_{E}}
\ar[d]^{\alpha^{0,d+1,1}_{(\mu_c(\overline{C}),\overline{D}),E}}
&&&\mu_{c+d}(\overline{C},\overline{D})\Box{E}
\ar[d]^{\alpha^{0,c+d,1}_{(\overline{C},\overline{D}),E}}\\
\mu_{d+2}(\mu_c(\overline{C}),\overline{D},E)
\ar[rrr]^{\alpha^{0,c,d+1}_{\overline{C},(\overline{D},E)}}
&&&\mu_{c+d+1}(\overline{C},\overline{D},E)
}$$
commutes. By (*) and the inductive definition of $\alpha^{0,i,j}$ this diagram can be replaced and expanded into the
following diagram
$$\xymatrix{
\mu_{2}(\mu_c(\overline{C}),\mu_d(\overline{D}))\Box E
\ar[rrr]^{\alpha^{0,c,1}_{\overline{C},\mu_d(\overline{D})}\Box id_{E}}
\ar[d]^{\alpha^{0,2,1}_{(\mu_c(\overline{C}),\mu_d(\overline{D})),E}}
&&&\mu_{c+1}(\overline{C},\mu_d(\overline{D}))\Box{E}
\ar[d]^{\alpha^{0,c+1,1}_{(\overline{C},\mu_d(\overline{D})),E}}\\
\mu_{3}(\mu_c(\overline{C}),\mu_d(\overline{D}),E)
\ar[rrr]^{\alpha^{0,c,2}_{\overline{C},(\mu_d(\overline{D}),E)}}
\ar[d]^{\mu_2(id_{\mu_c(\overline{C})},\alpha^{0,d,1}_{\overline{D},E})}
&&&\mu_{c+2}(\overline{C},\mu_d(\overline{D}),E)
\ar[d]^{\mu_{c+1}(id_{\overline{C}},\alpha^{0,d,1}_{\overline{D},E})}
\\
\mu_2(\mu_c(\overline{C}),\mu_{d+1}(\overline{D},E))
\ar[rrr]^{\alpha^{0,c,1}_{\overline{C},\mu_{d+1}(\overline{D},E)}}
&&&\mu_{c+1}(\overline{C},\mu_{d+1}(\overline{D},E))
}$$
The lower square commutes by naturality of $\alpha$. So it suffices to show the upper square commutes. This is just the
previous diagram with $\overline{D}$ replaced by $\mu_d(\overline{D})$. Thus we have reduced to the case $d=1$. We will
find it convenient to display this diagram in reflected form:
$$\xymatrix{
\mu_{2}(\mu_c(\overline{C}),D)\Box E
\ar[d]^{\alpha^{0,c,1}_{\overline{C},D}\Box id_{E}}
\ar[rrr]^{\alpha^{0,2,1}_{(\mu_c(\overline{C}),D),E}}
&&&\mu_{3}(\mu_c(\overline{C}),D,E)
\ar[d]^{\alpha^{0,c,2}_{\overline{C},(D,E)}}\\
\mu_{c+1}(\overline{C},D)\Box{E}
\ar[rrr]^{\alpha^{0,c+1,1}_{(\overline{C},D),E}}
&&&\mu_{c+2}(\overline{C},D,E)
}$$
We have to show this diagram commutes, where $D$ and $E$ are objects of $\scC$ and $\overline{C}$ is an object of $\scC^c$.

This diagram commutes trivially if $c\le 1$. So assume $c>1$ and $\overline{C}=(C_1,\overline{C}')$.
Again using (*) and the inductive definition of $\alpha^{0,i,j}$, we can expand this diagram into
$$\xymatrix{
\mu_{2}(\mu_2(C_1,\mu_{c-1}(\overline{C}')),D)\Box E
\ar[d]^{\alpha^{0,2,1}_{\mu_2(C_1,\mu_{c-1}(\overline{C}')),D}\Box id_{E}}
\ar[rrr]^{\alpha^{0,2,1}_{(\mu_2(C_1,\mu_{c-1}(\overline{C}'),D),E}}
&&&\mu_{3}(\mu_2(C_1,\mu_{c-1}(\overline{C}')),D,E)
\ar[d]^{\alpha^{0,2,2}_{(C_1,\mu_{c-1}(\overline{C}')),(D,E)}}\\
\mu_{3}(C_1,\mu_{c-1}(\overline{C}'),D)\Box{E}\ar@{=}[d]
\ar[rrr]^{\quad\alpha^{0,3,1}_{(C_1,\mu_{c-1}(\overline{C}'),D),E}}
&&&\mu_{4}(C_1,\mu_{c-1}(\overline{C}'),D,E)\ar@{=}[d]
\\
\mu_{2}(C_1,\mu_2(\mu_{c-1}(\overline{C}'),D))\Box{E}
\ar[d]^{\mu_2(id_{C_1},\alpha^{0,1,1}_{\mu_{c-1}(\overline{C}'),D})\Box id_E}
&&&\mu_{2}(C_1,\mu_3(\mu_{c-1}(\overline{C}'),D,E))
\ar[d]^{\mu_2(id_{C_1},\alpha^{0,1,2}_{\mu_{c-1}(\overline{C}'),(D,E)})}
\\
\mu_{2}(C_1,\mu_c(\overline{C}',D))\Box E
\ar@{=}[d]
&&&\mu_{2}(C_1,\mu_{c+1}(\overline{C}',D,E))
\ar@{=}[d]
\\
\mu_{c+1}(\overline{C},D)\Box{E}
\ar[rrr]^{\alpha^{0,c+1,1}_{(\overline{C},D),E}}
&&&\mu_{c+2}(\overline{C},D,E)
}$$
The top square in this diagram is the original diagram with $\overline{C}$ replaced by $(C_1,\mu_{c-1}(\overline{C}'))$, thus
reducing it to the case $c=2$. This top square can be expanded into the pentagonal diagram of the hypothesis of the theorem
and thus commutes. It remains to show that the bottom square commutes.

After rewriting the bottom square in reflected form and applying the inductive definition of $\alpha^{0,i,j}$ we obtain
the following expanded diagram
$$\xymatrix{
(C_1\Box(\mu_{c-1}(\overline{C}')\Box D))\Box{E}
\ar[rrr]^{\quad(id_{C_1}\Box\alpha^{0,1,1}_{\mu_{c-1}(\overline{C}'),D})\Box id_E}
\ar[d]^{\eta_{C_1,\mu_{c-1}(\overline{C}')\Box D,E}}
&&&(C_1\Box\mu_c(\overline{C}',D))\Box E
\ar[d]^{\eta_{C_1,\mu_c(\overline{C}',D),E}}
\\
C_1\Box((\mu_{c-1}(\overline{C}')\Box D)\Box E)
\ar[rrr]^{\quad id_{C_1}\Box(\alpha^{0,1,1}_{\mu_{c-1}(\overline{C}'),D}\Box id_E)}
\ar[d]^{id_{C_1}\Box\alpha^{0,2,1}_{(\mu_{c-1}(\overline{C}'),D),E}}
&&&C_1\Box(\mu_c(\overline{C}',D)\Box E)
\ar[d]^{id_{C_1}\Box\alpha^{0,c,1}_{(\overline{C}',D),E}}
\\
C_1\Box\mu_3(\mu_{c-1}(\overline{C}'),D,E)
\ar[rrr]^{\qquad id_{C_1}\Box\alpha^{0,c-1,2}_{\overline{C}',(D,E)}}
&&&C_1\Box\mu_{c+1}(\overline{C}',D,E)
}$$
The top square commutes by naturality of $\eta$ and the bottom square by naturality of $\alpha$.

This completes the verification that we have constructed an $A_\infty$-monoidal structure on $\scC$.
\end{proof}

To obtain the full version of LaPlaza's coherence theorem, we start with an operadic
reformulation of Theorem \ref{CohThm0}.

\begin{defi}
The LaPlaza operad $\scL=\{\scL_m\}_{m\ge0}$ is the operad in $CAT$ which acts on directed monoidal
categories as in the hypothesis of Theorem \ref{CohThm0}.  Specifically $\scL_m$ can be described as
a full subcategory of the free directed monoidal category on $m$ generating objects
$\{x_1, x_2,\dots, x_m\}$, whose objects look like $x_1\Box  x_2\Box\dots\Box x_m$ after removing
all parentheses. Thus $\scL_0=\{0\}$, $\scL_1=\{x_1\}$, and for $m\ge2$ the objects of $\scL_m$
are in bijective correspondence with planar binary trees with $m$ input edges.
\end{defi}

\begin{rema}
$\scL_2$ is the trivial poset $\{x_1\Box x_2\}$,
$\scL_3$ is the poset
$$\eta_{x_1,x_2,x_3}: (x_1\Box x_2)\Box x_3\longrightarrow x_1\Box(x_2\Box x_3),$$
isomorphic to $\scI$, while $\scL_4$ is the pentagonal poset generated by the labelled arrows shown below.
$$\xymatrix
{
&(x_1\Box x_2)\Box(x_3\Box x_4)\ar[dr]^{\eta_{x_1,x_2,x_3\Box x_4}}\\
((x_1\Box x_2)\Box x_3)\Box x_4
\ar[ur]^{\eta_{x_1\Box x_2,x_3,x_4}\quad}
\ar[d]_{\eta_{x_1,x_2,x_3}\Box id_{x_4}}
\ar[rr]
\ar[drr]|(.5){\quad\hole}
&& x_1\Box(x_2\Box(x_3\Box x_4))\\
(x_1\Box(x_2\Box x_3))\Box x_4\ar[rr]_{\eta_{x_1,x_2\Box x_3,x_4}}
\ar[urr]
&&x_1\Box((x_2\Box x_3)\Box x_4)\ar[u]_{\ id_{x_1}\Box\eta_{x_2,x_3,x_4}}
}$$
LaPlaza's coherence theorem states that $\scL_m$ is a poset for all $m$.
\end{rema}

\begin{rema} LaPlaza works with natural transformations $\eta_{A,B,C}:A\Box (B\Box C)\longrightarrow (A\Box B)\Box C$. Moreover,
he does not consider units. So, if $\scL^\ast \subset \scL$ is the suboperad obtained from $\scL$ by dropping the unit,
the operad $\scL^\ast$ is dual to LaPlaza's original one. Note also that $\scL_m=\scL^\ast_m$ for $m\geq 1$.
Our poset $\scL_m$, $m\ge 3$ is precisely the poset considered by Tamari \cite{Tam2}, and is now commonly
called the {\it Tamari lattice\/} \cite{Wiki}.
\end{rema}

With this notation, we can reformulate Theorem \ref{CohThm0} as follows.

\begin{theo}\label{CohThm1}
There is a map of $CAT$-operads
$$\Lambda: \scK\longrightarrow\scL$$
which is a surjection.
\end{theo}

The existence of $\Lambda$ is clear from the statement of Theorem \ref{CohThm0}.  Surjectivity
follows from the proof of Theorem \ref{CohThm0}, where it is shown that
$$(x_1x_2)x_3\longrightarrow x_1x_2x_3$$
maps via $\Lambda$ to
$$\eta_{x_1,x_2,x_3}: (x_1\Box x_2)\Box x_3\longrightarrow x_1\Box(x_2\Box x_3),$$
and the fact that $\eta_{x_1,x_2,x_3}$ generates $\scL$ as a $CAT$-operad.

\vspace{2ex}

LaPlaza's coherence theorem is not immediately apparent from Theorem \ref{CohThm1},
since a quotient category of a poset need not be a poset.  We need the following
additional observation.

\begin{lem}\label{inverse_images}
For any object $T\in\scL_m$, the inverse image under $\Lambda$ of the subcategory $\{T\}$
is a subposet of $\scK_m$ containing both a minimal and a maximal object.
\end{lem}

\begin{proof}
For $m=0,1,2$, the functor $\Lambda$ is an isomorphism and there is nothing to prove.  For $m\ge 3$, we
may regard $T$ as a planar binary tree. Clearly the minimal object of $\Lambda^{-1}\{T\}$  is
$T$ regarded as an object of $\scK_m$. The maximal object of $\Lambda^{-1}\{T\}$ is obtained
from $T$ by successively shrinking the rightmost incoming edge to every node of $T$, with the exception of
those edges which are leaves, till the rightmost edge of each node is a leaf.
\end{proof}

\begin{exam}
The inverse images in Lemma \ref{inverse_images} for $\Lambda:\scK_4\longrightarrow\scL_4$ are as follows:
\begin{eqnarray*}
\Lambda^{-1}\{((x_1\Box x_2)\Box x_3)\Box x_4\} &= &\left\{
((x_1x_2)x_3)x_4=
\raisebox{10pt}{$\xymatrix@=5pt@M=-1pt@W=-1pt{
\ar@{-}[dr] &&\ar@{-}[dl]&&&\ar@{-}[ddll] &&&\ar@{-}[dddlll]\\
&{\scriptscriptstyle\bullet}\ar@{-}[drr]\\
&&&{\scriptscriptstyle\bullet}\ar@{-}[drr]\\
&&&&&{\scriptscriptstyle\bullet}\ar@{-}[d]\\
&&&&&
}$}
\right\}\\
\Lambda^{-1}\{(x_1\Box x_2)\Box (x_3\Box x_4)\} &= &\left\{
\begin{picture}(250,20)
\put(0,0){$(x_1x_2)(x_3x_4)=$
\raisebox{10pt}{$\xymatrix@=5pt@M=-1pt@W=-1pt{
\ar@{-}[dr] &&\ar@{-}[dl]&&\ar@{-}[dr] &&\ar@{-}[dl]\\
&{\scriptscriptstyle\bullet}\ar@{-}[ddrr]&&&&{\scriptscriptstyle\bullet}\ar@{-}[ddll]\\
\\
&&&{\scriptscriptstyle\bullet}\ar@{-}[d]\\
&&&&&{}
}$}}
\put(120,4){\vector(1,0){20}}
\put(145,0){$(x_1x_2)x_3x_4=$
\raisebox{10pt}{$\xymatrix@=5pt@M=-1pt@W=-1pt{
\ar@{-}[dr]&&\ar@{-}[dl]&\ar@{-}[dddl]&&&\ar@{-}[dddllll]\\
&{\scriptscriptstyle\bullet}\ar@{-}[ddr]\\
\\
&&{\scriptscriptstyle\bullet}\ar@{-}[d]\\
&&
}$}}
\end{picture}
\right\}\\
\Lambda^{-1}\{(x_1\Box (x_2\Box x_3))\Box x_4\} &= &\left\{
\begin{picture}(250,20)
\put(0,0){$(x_1(x_2x_3))x_4=$
\raisebox{10pt}{$\xymatrix@=5pt@M=-1pt@W=-1pt{
\ar@{-}[ddrrr] &&&\ar@{-}[dr]&&\ar@{-}[dl] &&&\ar@{-}[dddllll]\\
&&&&{\scriptscriptstyle\bullet}\ar@{-}[dl]\\
&&&{\scriptscriptstyle\bullet}\ar@{-}[dr]\\
&&&&{\scriptscriptstyle\bullet}\ar@{-}[d]\\
&&&&
}$}}
\put(120,4){\vector(1,0){20}}
\put(145,0){$(x_1x_2x_3)x_4=$
\raisebox{10pt}{$\xymatrix@=5pt@M=-1pt@W=-1pt{
\ar@{-}[ddr]&\ar@{-}[dd]&\ar@{-}[ddl] &&&\ar@{-}[dddlll]\\
\\
&{\scriptscriptstyle\bullet}\ar@{-}[dr]\\
&&{\scriptscriptstyle\bullet}\ar@{-}[d]\\
&&
}$}}
\end{picture}
\right\}\\
\Lambda^{-1}\{x_1\Box ((x_2\Box x_3)\Box x_4)\} &= &\left\{
\begin{picture}(250,20)
\put(0,0){$x_1((x_2x_3)x_4)=$
\raisebox{10pt}{$\xymatrix@=5pt@M=-1pt@W=-1pt{
\ar@{-}[dddrrrr] &&&\ar@{-}[dr]&&\ar@{-}[dl] &&&\ar@{-}[ddlll]\\
&&&&{\scriptscriptstyle\bullet}\ar@{-}[dr]\\
&&&&&{\scriptscriptstyle\bullet}\ar@{-}[dl]\\
&&&&{\scriptscriptstyle\bullet}\ar@{-}[d]\\
&&&&
}$}}
\put(120,4){\vector(1,0){20}}
\put(145,0){$x_1(x_2x_3)x_4=$
\raisebox{10pt}{$\xymatrix@=5pt@M=-1pt@W=-1pt{
\ar@{-}[dddrrr]&&\ar@{-}[dr]&&\ar@{-}[dl]&&\ar@{-}[dddlll]\\
&&&{\scriptscriptstyle\bullet}\ar@{-}[dd]\\
\\
&&&{\scriptscriptstyle\bullet}\ar@{-}[d]\\
&&&
}$}}
\end{picture}
\right\}\\
\end{eqnarray*}
\begin{eqnarray*}
\Lambda^{-1}\{x_1\Box (x_2\Box (x_3\Box x_4))\} &=
&\begin{picture}(250,20)
\end{picture}
\end{eqnarray*}
$$\left\{
\begin{picture}(400,50)
\put(0,0){$x_1(x_2(x_3x_4))=$
\raisebox{10pt}{$\xymatrix@=5pt@M=-1pt@W=-1pt{
\ar@{-}[dddrrrr] &&&\ar@{-}[ddrr]&&\ar@{-}[dr] &&&\ar@{-}[dll]\\
&&&&&&{\scriptscriptstyle\bullet}\ar@{-}[dl]\\
&&&&&{\scriptscriptstyle\bullet}\ar@{-}[dl]\\
&&&&{\scriptscriptstyle\bullet}\ar@{-}[d]\\
&&&&
}$}}
\put(120,0){\vector(1,0){170}}
\put(120,5){\vector(2,1){40}}
\put(150,30){$x_1(x_2x_3x_4)=$
\raisebox{10pt}{$\xymatrix@=5pt@M=-1pt@W=-1pt{
\ar@{-}[dddrrr]&&&\ar@{-}[ddr]&\ar@{-}[dd]&\ar@{-}[ddl]\\
\\
&&&&{\scriptscriptstyle\bullet}\ar@{-}[dl]\\
&&&{\scriptscriptstyle\bullet}\ar@{-}[d]\\
&&&
}$}}
\put(120,-5){\vector(2,-1){40}}
\put(150,-35){$x_1x_2(x_3x_4)=$
\raisebox{10pt}{$\xymatrix@=5pt@M=-1pt@W=-1pt{
\ar@{-}[dddrrrr]&&&\ar@{-}[dddr]&&\ar@{-}[dr]&&\ar@{-}[dl]\\
&&&&&&{\scriptscriptstyle\bullet}\ar@{-}[ddll]\\
\\
&&&&{\scriptscriptstyle\bullet}\ar@{-}[d]\\
&&&&
}$}}
\put(250,25){\vector(2,-1){40}}
\put(250,-30){\vector(2,1){40}}
\put(300,0){$x_1x_2x_3x_4=$
\raisebox{10pt}{$\xymatrix@=5pt@M=-1pt@W=-1pt{
\ar@{-}[dddrrr]&&\ar@{-}[dddr]&&\ar@{-}[dddl]&&\ar@{-}[dddlll]\\
\\
\\
&&&{\scriptscriptstyle\bullet}\ar@{-}[d]\\
&&&
}$}}
\end{picture}
\right\}$$
\end{exam}

\begin{coro} (LaPlaza Coherence Theorem)
For all $m$, the category  $\scL_m$ is a poset  (known as the Tamari lattice for $m\ge3$).
\end{coro}

\begin{proof}
Let
$$\xymatrix{S\ar@<2pt>[rr]^f\ar@<-2pt>[rr]_g &&T}$$
be morphisms in $\scL_m$.  Since $\Lambda$ is surjective, we may find preimages
$$f':S'\longrightarrow T', \quad g':S''\longrightarrow T''$$
under $\Lambda$ of $f$ and $g$ respectively. Now let $S'''$ be the minimal element of
$\Lambda^{-1}\{S\}$ and let $T'''$ be the maximal element of $\Lambda^{-1}\{T\}$.
Then since $\scK_m$ is a poset, we have a commutative diagram
$$\xymatrix{&S'\ar[r]^{f'} &T'\ar[dr]\\
S'''\ar[ur]\ar[dr] &&&T'''\\
&S''\ar[r]^{g'} &T''\ar[ur]
}$$
Applying $\Lambda$ to this diagram and noting that $\Lambda$ sends the unlabelled arrows
to identities, we obtain $f=g$.
\end{proof}

We will now give an explicit description of the posets $\scL_m$ for $m\ge4$.  Similar considerations
apply to $\scK_m$ and give an alternative description of those posets.

\begin{defi}
Let $\scO$ be a $CAT$-operad with a single nullary operation $\scO_0=\{0\}$ (such as $\scK$ or $\scL$). Suppose $m\ge 4$
and let $\{a<b<c\}\subset\{1,2,3,\dots,m\}$. We define the functor 
$\pi_{a,b,c}:\scO_m\longrightarrow\scO_3$ to be the composite
$$\scO_m\longrightarrow\scO_m\times\prod_{i=1}^m\scO_{k_i}\longrightarrow\scO_3.$$
Here $k_i=0$ if $i \not\in\{a,b,c\}$, $k_a=k_b=k_c=1$, the first map takes $\phi\in\scO_m$ to
$(c;\epsilon_1,\epsilon_2,\dots,\epsilon_m)$, where $\epsilon_a=\epsilon_b=\epsilon_c=id\in\scO_1$
with all other $\epsilon_i=0\in\scO_0$, and the second map is composition in $\scO$.
\end{defi}

\begin{propo}
There is a commutative diagram
$$\xymatrix{
\scK_m\ar@{^(->}[rr]^(.3){\{\pi_{a,b,c}\}}\ar[d]^\Lambda &&\prod_{\{1\le a<b<c\le m\}}\scK_3\ar[d]^{\prod\Lambda}\\
\scL_m\ar@{^(->}[rr]^(.3){\{\pi_{a,b,c}\}} &&\prod_{\{1\le a<b<c\le m\}}\scL_3\ar[r]^\cong&\prod_{\{1\le a<b<c\le m\}}\scI
}$$
with the horizontal arrows being full imbeddings of posets.
\end{propo}

The proof is straight forward and left as an exercise for the reader.
\section{Rectification of $A_\infty$-monoidal categories}

It is well known that a monoidal category is equivalent to a strictly monoidal category, c.f. \cite[pages 257-259]{Mac2}.
[Recall that a monoidal category is strict if the associativity natural transformations $\eta_{A,B,C}$ of Theorem \ref{CohThm0}
are the identities.] We establish an analogous result for $A_\infty$-monoidal categories.

We first need a preliminary construction.

\begin{defi}
For $k\ge 2$ we define the poset $\hK_k$ to have as objects combinatorial trees as defined in \cite[Appendix E]{Lein} with $k$ input edges.
All nodes except the root node, i.e. the node at the output of the tree, are required to have more than one incoming edge. The root node 
may have zero, one, or more incoming edges.
 We define $T<T'$ if $T'$ can be obtained from $T$ by shrinking some internal edges. We define $\hK_1$ to consist of the single tree:
$$\xymatrix@=5pt@M=-1pt@W=-1pt{
\ar@{-}[dd]\\ \\
{\scriptscriptstyle\bullet}\ar@{-}[dd]\\ \\
\\
}$$
and $\hK_0$ to consist of the single tree
$$\xymatrix@=5pt@M=-1pt@W=-1pt{
{\scriptscriptstyle\bullet}\ar@{-}[ddd]\\ \\ \\
\\
}$$
\end{defi}
\begin{exam}
$$\xymatrix@=5pt@M=-1pt@W=-1pt{
\ar@{-}[dr] &&\ar@{-}[dl]\\
&{\scriptscriptstyle\bullet}\ar@{-}[d]\\
&{\scriptscriptstyle\bullet}\ar@{-}[d]\\
&&{}
}$$
is allowed in $\hK_2$, while
$$\xymatrix@=5pt@M=-1pt@W=-1pt{
&& \ar@{-}[d]\\
\ar@{-}[dr] && {\scriptscriptstyle\bullet}\ar@{-}[dl]\\
& {\scriptscriptstyle\bullet}\ar@{-}[d]\\
&&{}
}$$
is not allowed.
\end{exam}
The collection $\hK=\{\hK_k\}_{k\ge 0}$ is a right module over the associahedral operad, that is there are maps of posets
$$\hK_m\times\prod_{i=1}^m\scK_{k_i}\longrightarrow \hK_{k_1+k_2+\dots+k_m}$$
satisfying the usual associativity and unit conditions. This right action is defined in exactly the same way as we defined the operad
structure on $\scK$, with the single exception that when we compose with $0\in\scK_0$, we never delete the root node. 
Moreover, $\hK$ is also a left module over $Ass$, the trivial operad parametrizing strictly monoidal structures. The
left action
$$Ass(m)\times\prod_{i=1}^m\hK_{k_i}\cong \prod_{i=1}^m\hK_{k_i}\longrightarrow \hK_{k_1+k_2+\dots+k_m}$$
is given by\newline
\centerline{
\raisebox{-20pt}{$\left(\xymatrix@=5pt@M=-1pt@W=-1pt{\\ \\ \\ \\ \\ \\ \\}\right.$}
$\xymatrix@=5pt@M=-1pt@W=-1pt{
T_{i1} &T_{i2} &\dots &\ T_{ik_i}\\
\ar@{-}[ddrr]&\ar@{-}[ddr]& &\ar@{-}[ddl]\\
\\
&&{\scriptscriptstyle\bullet}\ar@{-}[dd]\\
\\
&&}$
\raisebox{-20pt}{$\left.\xymatrix@=5pt@M=-1pt@W=-1pt{\\ \\ \\ \\ \\ \\ \\}\right)_{i=1}^m\qquad\mapsto\qquad$}
$\xymatrix@=5pt@M=-1pt@W=-1pt{
T_{11} &\dots &T_{1k_1} &\dots &T_{m1} &\dots &T_{mk_m}\\
\ar@{-}[ddrrr]&&\ar@{-}[ddr]& &\ar@{-}[ddl]&&\ar@{-}[ddlll]\\
\\
&&&{\scriptscriptstyle\bullet}\ar@{-}[dd]\\
\\
&&&}$
}
It is clear that the left and right actions commute with each other, so $\hK$ is an $Ass$-$\scK$-bimodule.

\begin{theo}
There is a functorial construction $\scC\mapsto\scM\scC$ together with functors $I:\scC\longrightarrow\scM\scC$ 
and $E:\scM\scC\longrightarrow\scC$
which associates to each $A_\infty$-monoidal category $\scC$ a strictly monoidal category $\scM\scC$ such that
\begin{enumerate}
\item the induced maps by $I$ and $E$ on the nerves of the categories are mutually inverse homotopy equivalences, 
\item the functor $I$ induces a lax homomorphism of $A_\infty$-spaces in the sense of \cite{BV},
\item if $\scC$ is strictly monoidal then $E$ is a strictly monoidal functor.

\end{enumerate}
\end{theo}

\begin{proof}
Let
$$\scM\scC=\hK\otimes_{\scK}\scC=\left(\coprod_{k\ge 0}\hK_k\times\scC^k\right)/\approx$$
where the equivalence relation is given by
$$\left(T\circ(S_1,S_2,\dots,S_m),(\overline{A}_1,\overline{A}_2,\dots,\overline{A}_m\right)
\approx\left(T,S_1(\overline{A}_1),S_2(\overline{A}_2),\dots,S_m(\overline{A}_m)\right),$$
where $T\in\hK_m$, $S_i\in\scK_{k_i}$, $\overline{A}_i\in\scC^{k_i}$, for $i=1,2,\dots,m$.
The left action of $Ass$ on $\hK$ then induces a strict monoidal structure on $\scM\scC$.

There are functors $I:\scK\longrightarrow\hK$ and $J:\hK\longrightarrow\scK$. The
functor $I$ takes a tree $S\in\scK$ to the tree
$$\xymatrix@=5pt@M=-1pt@W=-1pt{
S\\
\ar@{-}[dd]\\ \\
{\scriptscriptstyle\bullet}\ar@{-}[dd]\\ \\
\\
}$$
The functor $E$ takes a tree in $\hK_k$, deletes the root vertex if it has only one or no incoming edges,
and regards it as a tree in $\scK$.  The composite $EI$ is the identity of $\scK$. There is a natural
transformation from the composite $IE$ to the identity of $\hK$, given by shrinking the edge above the
root vertex. The functors $I$, $E$ and this natural transformation are compatible with the right actions of
$\scK$ on $\hK$ and on itself.  Hence there are induced functors
$$I:\scC\longrightarrow\scM\scC, \qquad E:\scM\scC\longrightarrow\scC$$
such that $EI$ is the identity of $\scC$ and we have an induced natural transformation from $IE$ to
the identity of $\scM\scC$.  It follows that the maps induced by $I$ and $E$ on the nerves of these categories are mutually inverse
equivalences. 

Moreover the nerve of $\scM\scC$ is homeomorphic to the topological construction described
in \cite[Theorem 1.26]{BV} applied to the nerve of $\scC$. There it is shown that the map induced by $I$
is a lax homomorphism of $A_\infty$-spaces. The fact that $E$ is strictly monoidal, when $\scC$ is strictly
monoidal, is straight forward.
\end{proof}

%
%
%
%

\end{document}